\documentclass{amsproc}

\usepackage{dsfont,tikz,hyperref}

\newtheorem{theorem}{Theorem}[section]

\theoremstyle{definition}

\theoremstyle{remark}

\numberwithin{equation}{section}

\begin{document}

\title{Lessons I Learned from Richard Stanley}


\author{James Propp}
\address{UMass Lowell, 1 University Avenue, Lowell, MA 01854}
\curraddr{}
\email{}
\thanks{}
\dedicatory{to Richard Stanley, on the occasion of his 70th birthday}


\subjclass[2000]{Primary }

\date{}

\begin{abstract}
I will share with the reader what I have learned from Richard Stanley
and the ways in which he has contributed to research in combinatorics
conducted by me and my collaborators.
\end{abstract}

\maketitle


\bibliographystyle{amsplain}

\newcommand{\mR}{\mathds{R}}

\section{Two big ideas}

The biggest lesson I learned from Richard Stanley's work is,
{\it combinatorial objects want to be partially ordered!}
By which I mean: if you are trying to understand 
some class of combinatorial objects,
you should look at ways of putting a partial order on the class,
in hopes of finding one that has especially nice properties.
You won't always succeed, but when you do,
the gains are likely to more than justify the effort.

A related lesson that Stanley has taught me is,
{\it combinatorial objects want to belong to polytopes!}
That is: If you can find a way to view the objects you're interested in
as the vertices or facets of a polytope,
or as the faces (of all dimensions) of a polytope,
or as the lattice points inside a polytope,
then geometrical methods will give you a lot of combinatorial insight.

\section{Tilings and perfect matchings}

The two articles of Stanley's that had the greatest impact 
on my research were \cite{S86c} 
and \cite{S85}, which deal respectively with 
rhombus tilings of hexagons (Stanley calls them plane partitions 
whose three-dimensional diagram fits inside a box) and 
domino tilings of rectangles
(Stanley, taking the dual point of view, calls them dimer covers).

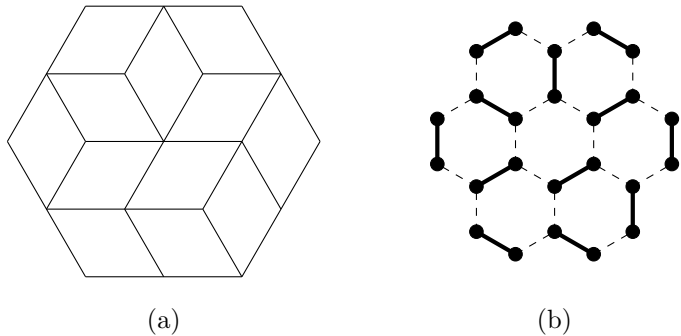
\begin{figure}[!h]
\centering
\begin{tikzpicture}[scale=0.3]
\path ({2*sqrt(3)},0) coordinate (A);
\path ({4*sqrt(3)},0) coordinate (B);
\path ({6*sqrt(3)},0) coordinate (C);
\path ({1*sqrt(3)},3) coordinate (D);
\path ({3*sqrt(3)},3) coordinate (E);
\path ({5*sqrt(3)},3) coordinate (F);
\path ({7*sqrt(3)},3) coordinate (G);
\path ({0*sqrt(3)},6) coordinate (H);
\path ({2*sqrt(3)},6) coordinate (I);
\path ({4*sqrt(3)},6) coordinate (J);
\path ({6*sqrt(3)},6) coordinate (K);
\path ({8*sqrt(3)},6) coordinate (L);
\path ({1*sqrt(3)},9) coordinate (M);
\path ({3*sqrt(3)},9) coordinate (N);
\path ({5*sqrt(3)},9) coordinate (O);
\path ({7*sqrt(3)},9) coordinate (P);
\path ({2*sqrt(3)},12) coordinate (Q);
\path ({4*sqrt(3)},12) coordinate (R);
\path ({6*sqrt(3)},12) coordinate (S);
\draw (A) -- (B) -- (C);
\draw (D) -- (E) -- (F);
\draw (I) -- (J) -- (K);
\draw (M) -- (N);
\draw (O) -- (P);
\draw (Q) -- (R) -- (S);
\draw (C) -- (G) -- (L);
\draw (F) -- (K) -- (P);
\draw (E) -- (J) -- (O);
\draw (D) -- (I);
\draw (N) -- (R);
\draw (H) -- (M) -- (Q);
\draw (A) -- (D) -- (H);
\draw (B) -- (E);
\draw (I) -- (M);
\draw (C) -- (F);
\draw (J) -- (N);
\draw (G) -- (K);
\draw (O) -- (R);
\draw (L) -- (P) -- (S);
\path ({13*sqrt(3)},1) coordinate (a);
\path ({15*sqrt(3)},1) coordinate (b);
\path ({12*sqrt(3)},2) coordinate (c);
\path ({14*sqrt(3)},2) coordinate (d);
\path ({16*sqrt(3)},2) coordinate (e);
\path ({12*sqrt(3)},4) coordinate (f);
\path ({14*sqrt(3)},4) coordinate (g);
\path ({16*sqrt(3)},4) coordinate (h);
\path ({11*sqrt(3)},5) coordinate (i);
\path ({13*sqrt(3)},5) coordinate (j);
\path ({15*sqrt(3)},5) coordinate (k);
\path ({17*sqrt(3)},5) coordinate (l);
\path ({11*sqrt(3)},7) coordinate (m);
\path ({13*sqrt(3)},7) coordinate (n);
\path ({15*sqrt(3)},7) coordinate (o);
\path ({17*sqrt(3)},7) coordinate (p);
\path ({12*sqrt(3)},8) coordinate (q);
\path ({14*sqrt(3)},8) coordinate (r);
\path ({16*sqrt(3)},8) coordinate (s);
\path ({12*sqrt(3)},10) coordinate (t);
\path ({14*sqrt(3)},10) coordinate (u);
\path ({16*sqrt(3)},10) coordinate (v);
\path ({13*sqrt(3)},11) coordinate (w);
\path ({15*sqrt(3)},11) coordinate (x);
\draw[fill=black] (a) circle (0.3cm);
\draw[fill=black] (b) circle (0.3cm);
\draw[fill=black] (c) circle (0.3cm);
\draw[fill=black] (d) circle (0.3cm);
\draw[fill=black] (e) circle (0.3cm);
\draw[fill=black] (f) circle (0.3cm);
\draw[fill=black] (g) circle (0.3cm);
\draw[fill=black] (h) circle (0.3cm);
\draw[fill=black] (i) circle (0.3cm);
\draw[fill=black] (j) circle (0.3cm);
\draw[fill=black] (k) circle (0.3cm);
\draw[fill=black] (l) circle (0.3cm);
\draw[fill=black] (m) circle (0.3cm);
\draw[fill=black] (n) circle (0.3cm);
\draw[fill=black] (o) circle (0.3cm);
\draw[fill=black] (p) circle (0.3cm);
\draw[fill=black] (q) circle (0.3cm);
\draw[fill=black] (r) circle (0.3cm);
\draw[fill=black] (s) circle (0.3cm);
\draw[fill=black] (t) circle (0.3cm);
\draw[fill=black] (u) circle (0.3cm);
\draw[fill=black] (v) circle (0.3cm);
\draw[fill=black] (w) circle (0.3cm);
\draw[fill=black] (x) circle (0.3cm);
\draw [ultra thick] (a) -- (c);
\draw [ultra thick] (b) -- (d);
\draw [ultra thick] (f) -- (j);
\draw [ultra thick] (g) -- (k);
\draw [ultra thick] (e) -- (h);
\draw [ultra thick] (i) -- (m);
\draw [ultra thick] (n) -- (q);
\draw [ultra thick] (l) -- (p);
\draw [ultra thick] (o) -- (s);
\draw [ultra thick] (r) -- (u);
\draw [ultra thick] (t) -- (w);
\draw [ultra thick] (v) -- (x);
\draw [dashed] (a) -- (d);
\draw [dashed] (b) -- (e);
\draw [dashed] (c) -- (f);
\draw [dashed] (d) -- (g);
\draw [dashed] (f) -- (i);
\draw [dashed] (g) -- (j);
\draw [dashed] (h) -- (k);
\draw [dashed] (h) -- (l);
\draw [dashed] (j) -- (n);
\draw [dashed] (k) -- (o);
\draw [dashed] (m) -- (q);
\draw [dashed] (n) -- (r);
\draw [dashed] (o) -- (r);
\draw [dashed] (p) -- (s);
\draw [dashed] (q) -- (t);
\draw [dashed] (s) -- (v);
\draw [dashed] (u) -- (w);
\draw [dashed] (u) -- (x);
\path ( {4*sqrt(3)},-2) node (y) {(a)};
\path ({14*sqrt(3)},-2) node (z) {(b)};
\end{tikzpicture}
\caption{A lozenge tiling and its dual perfect matching.}
\end{figure}
Figure 1(a) shows one of the 20 ways to tile
a regular hexagon of side-length 2
using twelve unit-rhombus tiles;
Figure 1(b) shows the associated perfect matching
of the graph whose edges correspond to
allowed positions of the tiles,
with vertices corresponding to triangular ``half-tiles''.

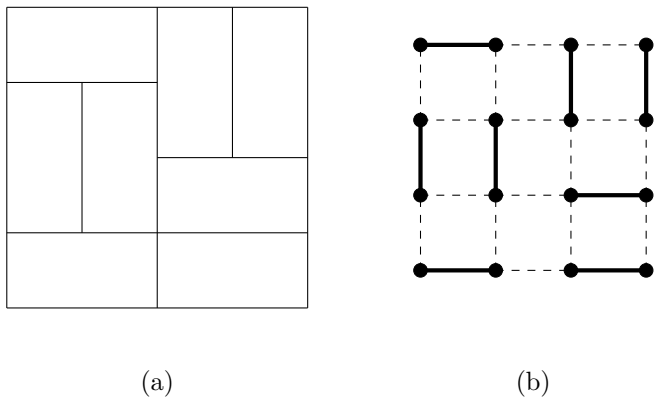
\begin{figure}[!h]
\centering
\begin{tikzpicture}[scale=0.5]
\path (0,0) coordinate (A);
\path (2,0) coordinate (B);
\path (4,0) coordinate (C);
\path (6,0) coordinate (D);
\path (8,0) coordinate (E);
\path (0,2) coordinate (F);
\path (2,2) coordinate (G);
\path (4,2) coordinate (H);
\path (6,2) coordinate (I);
\path (8,2) coordinate (J);
\path (0,4) coordinate (K);
\path (2,4) coordinate (L);
\path (4,4) coordinate (M);
\path (6,4) coordinate (N);
\path (8,4) coordinate (O);
\path (0,6) coordinate (P);
\path (2,6) coordinate (Q);
\path (4,6) coordinate (R);
\path (6,6) coordinate (S);
\path (8,6) coordinate (T);
\path (0,8) coordinate (U);
\path (2,8) coordinate (V);
\path (4,8) coordinate (W);
\path (6,8) coordinate (X);
\path (8,8) coordinate (Y);
\draw (A) -- (B) -- (C) -- (D) -- (E);
\draw (F) -- (G) -- (H) -- (I) -- (J);
\draw (M) -- (N) -- (O);
\draw (P) -- (Q) -- (R);
\draw (U) -- (V) -- (W) -- (X) -- (Y);
\draw (A) -- (F) -- (K) -- (P) -- (U);
\draw (G) -- (L) -- (Q);
\draw (C) -- (H) -- (M) -- (R) -- (W);
\draw (N) -- (S) -- (X);
\draw (E) -- (J) -- (O) -- (T) -- (Y);
\path (11,1) coordinate (a);
\path (13,1) coordinate (b);
\path (15,1) coordinate (c);
\path (17,1) coordinate (d);
\path (11,3) coordinate (e);
\path (13,3) coordinate (f);
\path (15,3) coordinate (g);
\path (17,3) coordinate (h);
\path (11,5) coordinate (i);
\path (13,5) coordinate (j);
\path (15,5) coordinate (k);
\path (17,5) coordinate (l);
\path (11,7) coordinate (m);
\path (13,7) coordinate (n);
\path (15,7) coordinate (o);
\path (17,7) coordinate (p);
\draw[fill=black] (a) circle (0.18cm);
\draw[fill=black] (b) circle (0.18cm);
\draw[fill=black] (c) circle (0.18cm);
\draw[fill=black] (d) circle (0.18cm);
\draw[fill=black] (e) circle (0.18cm);
\draw[fill=black] (f) circle (0.18cm);
\draw[fill=black] (g) circle (0.18cm);
\draw[fill=black] (h) circle (0.18cm);
\draw[fill=black] (i) circle (0.18cm);
\draw[fill=black] (j) circle (0.18cm);
\draw[fill=black] (k) circle (0.18cm);
\draw[fill=black] (l) circle (0.18cm);
\draw[fill=black] (m) circle (0.18cm);
\draw[fill=black] (n) circle (0.18cm);
\draw[fill=black] (o) circle (0.18cm);
\draw[fill=black] (p) circle (0.18cm);
\draw [ultra thick] (a) -- (b);
\draw [ultra thick] (c) -- (d);
\draw [ultra thick] (e) -- (i);
\draw [ultra thick] (f) -- (j);
\draw [ultra thick] (g) -- (h);
\draw [ultra thick] (k) -- (o);
\draw [ultra thick] (l) -- (p);
\draw [ultra thick] (m) -- (n);
\draw [dashed] (b) -- (c);
\draw [dashed] (e) -- (f);
\draw [dashed] (f) -- (g);
\draw [dashed] (i) -- (j);
\draw [dashed] (j) -- (k);
\draw [dashed] (k) -- (l);
\draw [dashed] (n) -- (o);
\draw [dashed] (o) -- (p);
\draw [dashed] (a) -- (e);
\draw [dashed] (b) -- (f);
\draw [dashed] (c) -- (g);
\draw [dashed] (d) -- (h);
\draw [dashed] (g) -- (k);
\draw [dashed] (h) -- (l);
\draw [dashed] (i) -- (m);
\draw [dashed] (j) -- (n);
\path ( 4,-2) node (y) {(a)};
\path (14,-2) node (z) {(b)};
\end{tikzpicture}
\caption{A domino tiling and its dual perfect matching.}
\end{figure}
Figure 2(a) shows one of the 36 ways to tile
a square of side-length 4
using eight 1-by-2 rectangular tiles (dominos);
Figure 2(b) shows the associated dimer cover (or perfect matching)
of the graph whose edges correspond to
allowed positions of the tiles,
with vertices corresponding to square half-tiles.

I was struck by the dissimilarity between 
formula (1) in \cite{S86c} and formula (2) in \cite{S85}.
The first of these implies that
the number of ways to tile a regular hexagon of side $n$
with $3n^2$ rhombuses of side length 1
(each consisting of two unit equilateral triangles joined edge-to-edge) is
$$\prod_{i=1}^{n} \prod_{j=1}^{n} \prod_{k=1}^{n} 
\frac{i+j+k-1}{i+j+k-2}\ .$$
The second formula implies that
the number of ways to tile a $2n$-by-$2n$ square 
with $2n^2$ dominos 
(each consisting of two unit squares joined edge-to-edge) is
$$\prod_{j=1}^{n} \prod_{k=1}^{n} 
\left( 4 \cos^2 \frac{\pi j}{2n+1} + 4 \cos^2 \frac{\pi k}{2n+1} \right).$$
Why do we see simple rational numbers in the former
and complicated trigonometric expressions in the latter,
when the two problems might seem at first to be so analogous to one another?
Pondering this question led me and others
into deeper exploration of the dimer model of statistical physics
(or what graph theorists call perfect matchings of graphs),
and M.I.T.\ became a major center for research in this field in the 1990s.
My sole coauthored paper with Stanley \cite{PS99}
was written during this period.
Trying to raise the visibility of this field, 
I took encouragement from \cite{S86b},
whose success emboldened me to come up with
a similar problems-list of my own \cite{P99}.
For an overview of the subject of enumeration of tilings, see \cite{P15}. 
For Stanley's own introduction to the subject, 
co-written with Federico Ardila
(a former member of the M.I.T.\ Tilings Research Group),
see \cite{AS10}.

A major tool in the study of tilings has been what are called height-functions.
These are mathematical constructions that generalize
a feature of tilings that you may have already noticed:
the human visual system is inclined to view Figure 1(a) as a projection 
of a stepped surface composed of squares seen at an oblique angle.
Your eye and brain may not perform the same trick for Figure 2(b),
but domino tilings are equally susceptible to being viewed
as surfaces in three-dimensions, from a purely mathematical perspective.
I learned about this point of view from
work of Conway and Lagarias \cite{CL90} and Thurston \cite{Th90},
and did some unpublished work showing how the key ideas
could be applied to a variety of combinatorial models \cite{P93}.
What is really going on is that the set of perfect matchings of a planar graph
can be endowed with a partial ordering
that turns it into a distributive lattice;
and, being at MIT, I was optimally situated to exploit this.
The distributive lattice structure on tilings
often gives the right way to ``$q$-ify'' enumerative questions.
My Ph.D.\ student David Wilson
came up with a brilliant way to exploit the lattice structure
to make it possible to efficiently sample from
the uniform distribution on the set of perfect matchings of a planar graph.
This is the method of Coupling From The Past (or CFTP),
originally developed for the study of tilings
but applicable much more broadly 
(for an overview, see \cite{PW98} or Chapter 22 of \cite{LPW09}).

Stanley's work on enumerating symmetry classes of plane partitions 
played an interesting role in the advent of the notion of cyclic sieving.
In \cite{S86c}, Stanley introduced the idea of complementing
a plane partition whose solid Young diagram fits inside a specific box,
and combined this new symmetry with other sorts of symmetry
that MacMahon and others had already studied.  John Stembridge,
in exploring some of Stanley's new symmetry classes, 
noticed a curious relation 
between these new enumeration problems and some old ones.
Specifically, he noticed that if one took
the $q$-enumeration of some old symmetry class, and set $q=-1$,
the result would be the number of plane partitions
which in addition to belonging to the symmetry class
also were self-complementary.
This is the $q=-1$ phenomenon of Stembridge \cite{S94}.
Vic Reiner, Dennis Stanton, and Dennis White went on to realize 
that the $q=-1$ phenomenon of Stembridge is just a special case 
of the more general cyclic sieving phenomenon \cite{RSW04}.
To be brief 
(though somewhat at variance with standard definitions and notation), 
we say that the quadruple $(S,\pi,p(x),\zeta)$
exhibits the cyclic sieving phenomenon
(with $S$ a finite set, $\pi$ a permutation of $S$,
$p(x)$ a polynomial with integer coefficients,
and $\zeta$ a root of unity) 
when for all integers $k$, the number of fixed points of $\pi^k$
equals $|p(\zeta^k)|$.  
For recent discussions of cyclic sieving, see \cite{RSW14} and \cite{S11}.

\section{Combinatorial reciprocity}

I have also been inspired by Stanley's work on combinatorial reciprocity,
as described in \cite{S73} and \cite{S74}.
(This aspect of Stanley's work was also the subject
of my presentation at the 2004 Stanley Conference \cite{P04}.)
The key result of the former paper
is that the chromatic polynomial of a graph $G$, evaluated at $-1$,
equals $(-1)^{|V(G)|}$ times the number of acyclic orientations of $G$.
When I first encountered this fact, it seemed miraculous.
In what sense do the acyclic orientation of $G$
correspond to colorings with $-1$ colors?
The answer is that one should think of
the set of colors as being the set of real numbers!
Schanuel \cite{S90} has taught us that 
the combinatorial Euler measure of $\mR$ is $-1$
(and more generally the combinatorial Euler measure 
of any polyhedral set homeomorphic to $\mR^k$ is $(-1)^k$),
so the preceding sentence is not entirely strange.
If one views an $\mR$-coloring of $G$ as a point in $\mR^{|V(G)|}$,
then the set of $\mR$-colorings of $G$ becomes
the complement of a hyperplane arrangement,
and we can view it as a disjoint union of $|V(G)|$-dimensional open cells
(each homeomorphic to $\mR^{|V(G)|}$),
which are in natural bijection with the acyclic orientations of $G$.
(For a related viewpoint, see \cite{ER98}.)

\cite{S74} describes many other examples in which
one starts with some polynomial $p(t)$
for which $p(n)$ has some enumerative significance
when $n$ is a positive integer,
and finds that the values of $p(n)$ when $n$ is a negative integer
possess (up to sign) some sort of enumerative significance as well,
reminiscent of but different from
the enumerative significance of $p(n)$ when $n$ is positive.
My favorite example comes from Ehrhart theory:
If $\Pi$ is a compact convex polytope, 
and $p(t)$ is its Ehrhart polynomial,
so that $p(n)$ is the number of lattice points
in the $n$th dilation of $\Pi$ (for all $n \neq 1$), 
then $p(-n)$ is the number of lattice points
in the {\em interior} of the $n$th dilation of $\Pi$.

One can even apply reciprocity to domino tilings.
For fixed $k$, the number of domino tilings
of a $k$-by-$n$ rectangle (call it $T_k(n)$) 
satisfies a linear recurrence relation in $n$
that allows us to extend it to all integers $n$, regardless of sign.
It turns out that, with this extended definition of $T_k$,
we have $T_k(-2-n) = \pm T_k(n)$.
For a combinatorial explanation of this, see \cite{P01}.
I am convinced that there is a lot of important work yet to be done
in the area of combinatorial reciprocity,
and I hope to see Stanley's articles serve as a foundation
for future progress.

\section{Dynamical algebraic combinatorics}

The articles of Stanley's that I've drawn nourishment from
most recently are \cite{S86a} and \cite{S09}.
These articles fit into a growing body of work
that one might call dynamical algebraic combinatorics
(a field that arguably includes within its purview
the cyclic sieving phenomenon described earlier).
The first article takes the combinatorial operation
that turns antichains of a poset into order ideals of a poset
and lifts it into a piecewise-linear map between
the order polytope (whose vertices correspond to order ideals)
and the chain polytope (whose vertices correspond to antichains).
The second article treats, among other things,
the operation of promotion on linear extensions of a poset.
Here I will mention a link between the two articles,
discussed in greater length in \cite{EP14}.
Sch{\"u}tzenberger's promotion operator 
on the set of semistandard Young tableaux
of rectangular shape with $A$ rows and $B$ columns
having entries between 1 and $n$ is naturally conjugate to an action 
on the rational points in the order polytope of $[A] \times [n-A]$
with denominator dividing $B$.
(Here $[n]$ denotes the chain of length $[n]$.)
The latter action, introduced in \cite{SW12},
is expressible as a composition of fundamental involutions called toggles,
which in the setting of \cite{EP14} can be seen as 
continuous piecewise-linear maps from the order polytope to itself.

The notion of the order polytope has caught on, 
but the allied notion of the chain polytope has languished in comparison:
the two search terms generated 270 and 46 hits respectively
in {\tt scholar.google.com} in June 2014.
I hope the latter notion will attract more of the attention it deserves.

\section{Enumerative Combinatorics, volumes 1 and 2}

No discussion of Stanley's contributions would be complete
without mention of his books \cite{S12} and \cite{S99}.  These books
are not light reading, but they are clearly written and
loaded with useful information.  A good deal of my email
correspondence with Stanley over the past two decades
consists of me asking him a question and him informing
me that the answer to my question (or some new result of mine) 
is in \cite{S12}). 
I would estimate that over the course of my career thus far 
I've spent several dozen hours rediscovering things 
that were already in these books.

In a light-hearted vein, I expressed my appreciation for these books
in the form of a song that was performed at the opening day banquet
of the Stanley@70 conference in 2014.
It's based on the song ``Guys and Dolls'' by Frank Loesser,
and some of the lines were written by Noam Elkies,
who also did the arrangements and conducted the performance from the piano.
\begin{verse}
What's in {\it Inventiones}?  
I'll tell you what's in {\it Inventiones}. \\
Folks provin' theorems, 'stead-a figurin' odds 
to use for bettin' on the ponies. \\
That's what in {\it Inventionies\/}! \\
What's in the {\it Intelligencer}?  
I'll tell you what's in the {\it Intelligencer}. \\
Articles on abstruse mathematical questions 
for which countin' plays a role in the enswer. \\
That's what's in the {\it Intelligencer} ! \\
What's in every math joynal?  I'll tell you what's in every math joynal. \\
Combinatorics achievin' renown as a fountain of truths 
both beautiful and etoynal. \\
That's what's in every math joynal ! \\
Combinatorialists have one trusted resource; \\
And now it's both a physical and an ``e-''source. \\
Yes sir!  Yes ma'am! \\
{\ } \\
When you study balls stuck in separate stalls, \\
Then the facts that you need are in EC One. \\
When you seek a combinatorial truth, \\
EC One's where you go to see if it's so -- unless it's in Knuth. \\
When you see a mu with a zeta or two \\
And a delta thrown in for some extra fun ... \\
Call it odd, call it even; it's a principle to believe in \\
That the source you're consultin' is EC One. \\
{\ } \\
When you see a rook that determines a hook \\
Then the book that you're lookin' at's EC One. \\
When a theorem features a bent letter S \\
All curled up in distress, the theorem's address is not hard to guess. \\
When you wend a path and some elegant math \\
Tells the number of ways that it can be done, \\
It's a true proposition known to every mathematician \\
That the opus you've opened is EC One. \\
{\ } \\
When your lovely proof springs a leak in its roof \\
Then the patch for your goof is in EC One. \\
But proceed with caution.  You never can tell: \\
Maybe page eighty-six has not just the fix but {\it your} proof as well! \\
When an exercise makes you feel not so wise \\
'Cause for you it ain't ``EC'' -- forgive the pun -- \\
Call it plus, call it minus; chalk it all up to Stanley's slyness \\
'Cause the book that has stumped you is EC One -- \\
Or EC Two -- \\
The book that you're readin' is EC One!
\end{verse}

\section{Special sequences}

Of particular note is Stanley's impressive list of combinatorial
incarnations of the Catalan numbers \cite{S13}.
No other integer sequence, not even the Fibonacci sequence,
has such a rich assortment of seemingly unrelated manifestations,
and the On-line Encyclopedia of Integer Sequences page for the Catalan sequence
(entry A000108) is the longest in the whole OEIS database.
Given the length of Stanley's annotated list, the most time-efficient way 
for a researcher to find out whether some particular Catalan-incarnation
has been noticed before is not to leaf through the whole addendum
(96 pages at present count) but to ask Stanley. 
Unfortunately, that will not be a viable method forever.
As Sara Billey and Bridget Tenner have pointed out,
we might start trying to think now about how
to automate what Richard does when he fields a question
about the Catalan number literature
by finding ways to assign ``fingerprints''
to different sorts of combinatorial objects,
in a way that might make some form of automated search possible.
This problem isn't just of interest to combinatorialists;
one can argue that it is a natural test-bed
for the much broader enterprise of semantic search.
In identifying distinguishing structural characteristics
of different incarnations of the Catalan objects, 
and finding ways to represent these distinctions in software,
we will learn lessons that can be applied much more broadly
to the Artificial Intelligence problem of content-based searching
in other mathematical domains.

Another integer sequence that has followed Richard around
over the past decades
(though not as doggedly as the Catalan sequence)
is the sequence $1,2,8,64,1024,\dots$.
In \cite{S89}, Stanley (citing earlier work of Mills, Robbins, and Rumsey)
mentions that $2^{{n \choose 2}}$ occurs as the sum of $2^{s(T)}$,
where $T$ ranges over all $n$-by-$n$ alternating sign matrices
and $s(T)$ is the number of $-1$'s in $T$.
(For a combinatorial explanation of this, 
relating the formula to domino tilings, see \cite{EKLP92}.)
Quite recently, this sequence resurfaced in \cite{LS13}
in connection with a seemingly quite different sort of problem
arising from Ramsey theory.  Having absorbed the lesson
``combinatorial objects want to be partially ordered''
from Stanley during my many years in the Boston area,
it was gratifying to have a chance to return the favor;
my contribution to \cite{LS13} was suggesting the partial ordering
that was one of the keys to unlocking the problem.

\section{The culture of combinatorics}

Finally, I'd like to mention a contribution that Richard
has made to the mathematical life of the Boston area
that might not otherwise be recorded, namely, his founding
and continued leadership of the Cambridge Combinatorics Coffee Club.  
``CCCC'' (as it is called) has served as a great way
for Boston-area combinatorialists to keep up
with one another's research interests,
and a way for some of us to incubate our ideas over time.
I wonder whether other branches of mathematics
have similar institutions
in which they freely share work-in-progress
and welcome others to join their projects.
I suspect that one reason many of us who work in combinatorics
have gravitated toward the field is
how friendly and uncompetitive its practitioners are;
one rarely hears about combinatorialists 
racing one another towards the solution of some hot open problem
or stealing ideas from each other.

No doubt the accessibility of the subject matter of combinatorics
is a large factor in the friendliness of practitioners.
But a big part is played by people like Stanley,
who lead by example and set the tone for the field,
creating what Margaret Bayer aptly described 
(in her remarks at the close of the Stanley@70 conference)
as ``a culture of cooperation and openness''.

\bibliography{propp}

\end{document}